\documentclass[a4paper,11pt]{article}
\usepackage{latexsym}
\usepackage{amsmath}
\usepackage{amssymb}
\usepackage{amsfonts}
\usepackage[dvips]{graphicx}

\newcommand{\refset}{\noindent \hangindent=10pt \hangafter=1}
\newcommand{\btheta}{{\bar{\theta}}}
\newcommand{\dd}{\mbox{d}}
\newcommand{\ep}{\varepsilon}
\newcommand{\qed}{\hfill $\Box$}
\begin{document}
\begin{center}
{\Large Bayesian Predictive Densities} \\
{\Large Based on Latent Information Priors} \\

\vspace{20pt}
{\large Fumiyasu Komaki}

{\small Department of Mathematical Informatics} \\
{\small Graduate School of Information Science and Technology, the University of Tokyo} \\
{\small 7-3-1 Hongo, Bunkyo-ku, Tokyo 113-8656, JAPAN} \\
\end{center}

\vspace{10pt}
\begin{center}
{\bf Summary}
\end{center}

Construction methods for prior densities are investigated from a predictive viewpoint.
Predictive densities for future observables are constructed by using observed data.
The simultaneous distribution of future observables and observed data is assumed to
belong to a parametric submodel of a multinomial model.
Future observables and data are possibly dependent.
The discrepancy of a predictive density to the true conditional density of future observables
given observed data is evaluated by the Kullback-Leibler divergence.
It is proved that limits of Bayesian predictive densities form an essentially complete class.
Latent information priors are defined as priors
maximizing the conditional mutual information between the parameter and the future observables
given the observed data.
Minimax predictive densities are constructed as limits of Bayesian predictive densities
based on prior sequences converging to the latent information priors.

\vspace{0.5cm}


\vspace{0.3cm}
\refset
\section{\normalsize Introduction}

We construct predictive densities for future observables by using observed data.
Future observables and data are possibly dependent and the simultaneous distribution of them
is assumed to belong to a submodel of a multinomial model.
Various practically important models such as categorical models
and graphical models are included in this class.

Let ${\cal X}$ and ${\cal Y}$ be finite sets composed of $k$ and $l$ elements,
and let $x$ and $y$ be random variables that take values in  ${\cal X}$ and ${\cal Y}$, respectively.
Let ${\cal M} = \{p(x,y|\theta) \mid \theta \in \Theta\}$ be a set of probability densities on ${\cal X} \times {\cal Y}$.
The model ${\cal M}$ is regarded as a submodel of the $kl$-nominal model with trial number 1.
Here, we do not lose generality by assuming the trial number is $1$.
The model $\mathcal{M}$ is naturally regarded as a subset of the hyperplane 
$\{p=(p_{ij}) \mid \sum\limits_{i=1}^{k}  \sum\limits_{j=1}^{l} p_{ij} = 1 \}$ in Euclidean space $\mathbb{R}^{kl}$.
In the following, we identify $\Theta$ with ${\cal M}$.
Then, the parameter space $\Theta$ is endowed with the induced topology as a subset of $\mathbb{R}^{kl-1}$.

A predictive density $q(y;x)$ is defined as a function from $\mathcal{X} \times \mathcal{Y}$ to $[0,1]$
satisfying $\sum\limits_{y \in {\cal  Y}}q (y;x) = 1 ~(x \in {\cal X})$.
The closeness of $q(y;x)$ to the true conditional probability density $p(y|x,\theta)$ is evaluated by the average Kullback-Leibler divergence: 
\begin{align}
R(\theta,q) = \sum_{x,y} p(x,y|\theta) \log \frac{p(y|x,\theta)}{q(y;x)}, \label{R}
\end{align}
where we define
$c \log 0 = - \infty$ $(c > 0)$, $0 \log 0 = 0$, $0 \log (c/0) = 0$ $(c \geq 0)$.
Although the conditional probability $p(y|x,\theta)$ is not uniquely defined when $p(x|\theta)=0$,
the risk value $R(\theta,q)$ is uniquely determined
because $p(x,y|\theta) \log p(y|x,\theta) = 0$ if $p(x|\theta)=0$.

First, we show that, for every predictive density $q(y;x)$, there exists
a limit $\lim\limits_{n \to \infty} p_{\pi_n} (y;x)$ of Bayesian predictive densities
\[
p_{\pi_n}(y|x) := \frac{\int p(x,y|\theta) \dd \pi_n(\theta)}{\int p(x|\theta) \dd \pi_n(\theta)},
\]
where $\{\pi_n\}^\infty_{n=1}$ is a prior sequence,
such that $R (\theta, \lim\limits_{n \to \infty} p_{\pi_n} (y;x)) \le R (\theta, q(y;x))$
for every $\theta \in \Theta $.
In the terminology of statistical decision theory,
this means that the class of predictive densities that are limits of Bayesian predictive densities is an essentially complete class.

Next, we investigate latent information priors defined as priors maximizing
the conditional mutual information between $y$ and $\theta$ given $x$.
We obtain a constructing method for a prior sequence $\{\pi_n\}^\infty_{n=1}$ converging the latent information prior,
based on which a minimax predictive density $\lim\limits_{n \to \infty} p_{\pi_n} (y|x)$ is obtained.
We consider limits of Bayesian predictive densities to deal with conditional probabilities.

There exist important previous studies on prior construction
by using the unconditional mutual information.
The reference prior by Bernardo (1979), (2005)
is a prior maximizing the mutual information between $\theta$ and $y$
 in the limit of the amount of information of $y$ goes to infinity.
It corresponds to the Jeffreys prior if there are no nuisance parameters;
see Ibragimov and Hasminskii (1973) and Clarke and Barron (1994) for rigorous treatments.
In coding theory, the prior maximizing the mutual information between $y$ and $\theta$ is used
for Bayes coding.
It was shown that the Bayes codes for finite alphabet models based on the priors are minimax
by Gallager (1979) and Davisson and Leon-Garcia (1980).
In our framework, these settings correspond to prediction of $y$ without $x$.
In statistical applications, $x$ plays an important role because it corresponds to observed data,
although ${\cal X}$ is an empty set in the reference analysis and
the standard framework of information theory;
see also Komaki (2004) for the relation between statistical prediction and Bayes coding.

Geisser (1978), in the discussion of Bernardo (1978),
discussed minimax prediction based on the risk function \eqref{R} as an alternative to
the reference prior approach.

The latent information priors introduced in the present paper bridge these two approaches.
The theorems obtained below clarify the relation between the conditional mutual information and minimax prediction based on observed data.

For Bayesian prediction of future observables by using observed data,
Akaike (1983) discussed priors maximizing the mutual information
between $x$ and $y$ and called them minimum information priors.
Kuboki (1998) also proposed priors for Bayesian prediction based on an information theoretic quantity.
These priors are different from latent information priors investigated in the present paper.

In section 2, we prove that, for every predictive density $q(y;x)$,
there exists a predictive density that is a limit of Bayesian predictive densities whose performance is not worse than that of $q(y;x)$.
In section 3, we introduce a construction method for minimax predictive densities as limits of Bayesian predictive densities.
The method is based on the conditional mutual information between $y$ and $\theta$ given $x$.
In section 4, we give some numerical results and discussions.

\section{\normalsize Limits of Bayesian predictive densities}

In this section, we prove
that the class of predictive densities that are limits of Bayesian predictive densities is an essentially complete class.

Throughout this paper, we assume the following conditions:

\vspace{0.5cm}

\noindent
Assumption 1. $\Theta$ is compact.

\vspace{0.5cm}
\noindent
Assumption 2.
For every $ x \in \mathcal{X}$, there exists $\theta \in \Theta$ such that $p(x|\theta) > 0$.

\vspace{0.5cm}
\noindent
These assumptions are not restrictive.
For Assumption 1, if $\Theta$ is not compact,
we can regard the closure $\bar{\Theta}$ as the parameter space instead of $\Theta$
because we consider a submodel of a multinomial model.
We do not lose generality by Assumption 2 because we can adopt $\mathcal{X} \setminus \{{x}_0\}$
instead of $\mathcal{X}$ if there exists ${x}_0 \in \mathcal{X}$ such that $p({x}_0|\theta) = 0$ for every $\theta \in \Theta$.

We prepare several preliminary results to prove Theorem 1 below.

Let ${\cal P}$ be the set of all probability measures on $\Theta$ endowed with the weak convergence topology
and the corresponding Borel algebra.
By the Prohorov theorem and Assumption 1, ${\cal P}$ is compact.

When $x$ and $y$ are fixed, the function
$\theta \in \Theta \longmapsto  p(x,y|\theta) \in [0,1] $
is bounded and continuous.
Thus, for every fixed $(x,y) \in \mathcal{X} \times \mathcal{Y}$, the function
\[
\pi \in {\cal P}  \longmapsto  \displaystyle p_\pi(x,y) := \int p(x,y|\theta) \dd\pi( \theta)
\]
is continuous, because of the definition of weak convergence.
Therefore, for every predictive density $q(y;x)$, the function from $\mathcal{P}$ to $[0,\infty]$ defined by
\begin{align}
D_q (\pi) := &\sum_{x,y} p_\pi (x,y) \log \frac{p_\pi(x,y)}{q(y;x)p_\pi(x)} \notag \\
= &\sum_{x,y} p_\pi (x,y) \log p_\pi(x,y) - \sum_x p_\pi(x) \log p_\pi(x)
- \hspace{-10pt} \sum_{(x,y):q(y;x) > 0} \hspace{-10pt} p_\pi (x,y) \log q(y;x) \notag \\
& - \hspace{-10pt} \sum_{(x,y):q(y;x) = 0} \hspace{-10pt} p_\pi (x,y) \log q(y;x) \label{p}
\end{align}
is lower semicontinuous, because the last term in \eqref{p} is lower semicontinuous and the other terms are continuous.

\vspace{0.5cm}

\noindent
Lemma 1. 
Let $\mu$ be a probability measure on $\Theta$.
Then,
$\mathcal{P}_{\varepsilon \mu}  = \{ \varepsilon \mu + \left(1-\varepsilon \right) \pi \mid \pi \in \mathcal{P} \}$ $( 0 \le \varepsilon \le 1 )$
is a closed subset of $\mathcal{P}$.

\vspace{0.5cm}

\noindent
Proof.
Suppose that $\pi_\infty \in \mathcal{P} $ is the limit of a convergent sequence $\{\pi_k\}^\infty_{k=1}$ in $\mathcal{P}_{\varepsilon \mu}$.
Since $\pi_k \in \mathcal{P}_{\varepsilon \mu} $, 
\[
\int f(\theta) \dd\pi_{k} (\theta) - \varepsilon \int f (\theta) \dd \mu(\theta) \ge 0
\]
for every nonnegative bounded continuous function $f(\theta)$ on $\Theta$.
Thus, 
\[
\int f (\theta) \dd \pi_\infty(\theta) = \lim_{k \to \infty} \int f(\theta) \dd\pi_{k} (\theta) \ge  \varepsilon \int f(\theta) \dd\mu(\theta).
\]
Hence, $\pi_\infty - \varepsilon \mu $ is a nonnegative measure.
Therefore, $ \pi_\infty \in \mathcal{P}_{\varepsilon \mu} $, and  $ \mathcal{P}_{\varepsilon \mu}$ is a closed set in $\mathcal{P}$.
\qed

\vspace{0.5cm}

\noindent
Lemma 2.
Let $f(\cdot)$ be a continuous function from $\mathcal{P}$ to $[0,\infty]$,
and let $\mu$ be a probability measure on $\Theta$ such that $p_\mu (x) := \int p(x|\theta) \dd \mu(\theta) > 0$ for every $x \in \mathcal{X}$.
Then, there is a probability measure $\pi_n$ in
\[
\mathcal{P}_{\mu/n} := \left\{ \frac{1}{n} \mu + \left(1-\frac{1}{n}\right) \pi \bigg| \pi \in \mathcal{P} \right\} ~~~~~ (n=1,2,3,\ldots)
\]
such that $f (\pi_n) = \inf\limits_{\pi \in \mathcal{P}_{\mu/n}} f(\pi)$.
Furthermore, there exists a convergent subsequence $ \{ \pi'_m \}^\infty_{m=1}$
of $\{ \pi_{n} \}^\infty_{n=1}$
such that the equality $f(\pi'_\infty) = \inf\limits_{\pi \in \mathcal{P} } f(\pi)$ holds,
where $\displaystyle \pi'_\infty = \lim\limits_{m \to \infty} \pi'_m$.

\vspace{0.5cm}

\noindent
Proof.
Note that there exists $\mu \in \mathcal{P}$ such that $p_\mu (x) := \int p(x|\theta) \dd \mu(\theta) > 0$ for every $x \in \mathcal{X}$ by Assumption 2.
By Lemma 1,
the sets $\mathcal{P}_{\mu/n} $ $(n = 1,2,3, \dotsc)$ are compact because they are closed subsets of a compact set $\mathcal{P}$.
Thus, there is a probability measure $\pi_n$ in
$\mathcal{P}_{\mu/n}$ such that $f (\pi_n) = \inf\limits_{\pi \in \mathcal{P}_{\mu/n}} f(\pi)$.
There exists a convergent subsequence $\{ \pi'_m \} ^{\infty}_{m=1}$ of $\{ \pi_n \} ^{\infty}_{n=1}$
because $ \mathcal{P}$ is compact.

Since $\mathcal{P}$ is compact and $f(\pi)$ is a continuous function of $\pi \in \mathcal{P}$,
there exists $\hat{\pi} \in \mathcal{P}$
such that $f(\hat{\pi}) = \inf\limits_{\pi \in \mathcal{P}} f(\pi)$.
Thus, $f(\pi'_\infty) \ge f(\hat{\pi})$,
where $\pi'_\infty := \lim\limits_{m \rightarrow \infty} \pi'_m$.
For every $ \varepsilon > 0$, there exists $\delta > 0 $ such that  
$\sup\limits_{d(\hat{\pi}, \pi)< \delta} f(\pi) \le f(\hat{\pi}) + \varepsilon$,
where $d$ is the Prohorov metric on $\mathcal{P}$.
We put
\[
\hat{\pi}_n = \frac{1}{n} \mu + \frac{n-1}{n} \hat{\pi} ~~~~~ (n=1,2,3,\ldots).
\]
Then, $\hat{\pi}_n \in \mathcal{P}_{\mu/n}$ and $\lim\limits_{n \rightarrow \infty} \hat{\pi}_n = \hat{\pi}$.
Thus, for every $\delta > 0$, there exists a positive integer $N$ such that $d(\hat{\pi},\hat{\pi}_n) < \delta$ $(n \geq N)$.
If $n \ge N$, then $f(\pi'_{\infty}) \le  f(\pi_n) \le f(\hat{\pi}_n) \le f(\hat{\pi})+ \varepsilon $.
Since $\varepsilon > 0$ is arbitrary, we have $f(\pi'_\infty ) \leq f(\hat{\pi})$.
Therefore, $f(\pi'_\infty ) = f(\hat{\pi})$. \qed

\vspace{0.5cm}

The conditional probability $p_\pi(y|x)$ is not uniquely specified if $p_\pi(x) = 0$.
To resolve the problem, we consider a sequence of priors $\{\pi_n\}_{n=1}^\infty$
that satisfies $p_{\pi_n}(x) > 0$ for every $n$ and $x \in \mathcal{X}$.
In the following, $\lim\limits_{n \to \infty} p_{\pi_n} (y|x)$ is defined to be a map
from $(x,y) \in \mathcal{X} \times \mathcal{Y}$ to the limit of
the real number sequence $\{p_{\pi_n} (y|x)\}_{n=1}^\infty$.
If there exist limits of sequence of real numbers $\{p_{\pi_n} (y|x)\}_{n=1}^\infty$
for all $(x,y) \in \mathcal{X} \times \mathcal{Y}$,
we say the limit $\lim\limits_{n \to \infty} p_{\pi_n} (y|x)$ of Bayesian predictive densities exists.
Obviously, if the limit $\lim\limits_{n \to \infty} p_{\pi_n} (y|x)$ exists, it is a predictive density because
$0 \leq \lim\limits_{n \to \infty} p_{\pi_n} (y|x) \leq 1$ for every $(x,y) \in \mathcal{X} \times \mathcal{Y}$
and $\sum\limits_{y \in \mathcal{Y}} \lim\limits_{n \to \infty} p_{\pi_n} (y|x) = 1$ for every $x \in \mathcal{X}$.

\vspace{0.5cm}
\noindent  
Theorem 1.
\begin{itemize}
\item[1)]
For every predictive density $q(y;x) $,
there exists a convergent prior sequence $\{\pi_n\}^\infty_{n=1}$
such that the limit $\lim\limits_{n \to \infty} p_{\pi_n} (y|x) $ exists and 
$R (\theta, \lim\limits_{n \to \infty} p_{\pi_n} (y|x)) \le R(\theta, q(y;x)) $ for every $\theta \in \Theta $.
\item[2)]
If there exists $\hat{\pi} \in \mathcal{P}$ such that $D_q (\hat{\pi}) = \inf\limits_{\pi \in \mathcal{P}} D_q (\pi)$
and $ p_{\hat{\pi}} (x) > 0$ for every $x \in \mathcal{X}$,
then $R (\theta, p_{\hat{\pi}} (y|x)) \le R(\theta, q(y;x)) $ for every predictive density $q(y;x)$ and $\theta \in \Theta $.
\end{itemize}

\vspace{0.5cm}
\noindent
Proof. 1)~
Let $\mathcal{N}^q := \{ (x,y) \in \mathcal{X} \times \mathcal{Y} \mid q(y;x) = 0 \}$
and $\Theta^q := \{ \theta \in \Theta \mid \sum\limits_{(x,y) \in \mathcal{N}^q} p(x,y|\theta) = 0\}$.
Let $\mathcal{P}^q$ be the set of all probability measures on $\Theta^q$.
Then, $\Theta^q$ and $\mathcal{P}^q$ are compact subsets of $\Theta$ and $\mathcal{P}$, respectively.

If $\Theta^q = \emptyset$, the assertion is obvious,
because $R(\theta, q(y;x)) = \infty$ for $\theta \notin \Theta^q$.
We assume that $\Theta^q \neq \emptyset$ in the following.
Let $\mathcal{X}^q := \{ x \in \mathcal{X} \mid \exists \theta \in \Theta^q$ such that $P(x|\theta)> 0\}$
and $\mu^q$ be a probability measure on $\Theta^q$
such that $p_{\mu^q} (x) := \int p(x|\theta) \dd \mu^q(\theta) > 0$ for every $x \in \mathcal{X}^q $.

Then, because $D_q (\pi) $ defined by \eqref{p} as a function of $\pi \in \mathcal{P}^q$ is continuous,
there exists $\pi_n \in \mathcal{P}_{\mu^q/n}^q := \{(1/n) \mu^q + (1-1/n) \pi \mid \pi \in \mathcal{P}^q \}$
such that $D_q(\pi_n) = \inf\limits_{\pi \in \mathcal{P}_{\mu/n}^q} D_q (\pi)$.
From Lemma 2, there exists a convergent subsequence
$\{ \pi'_m \} ^{\infty}_{m=1}$ of $\{ \pi_{n} \} ^{\infty}_{n=1}$
such that $D_q (\pi'_\infty) =  \inf\limits_{\pi \in \mathcal{P}^q} D_q (\pi)$,
where $\pi'_\infty = \lim\limits_{m \rightarrow \infty} \pi'_m$.

Let $n_m$ be the integer satisfying $\pi'_m = \pi_{n_m}$.
We can take a subsequence $\{\pi'_m\}_{m=1}^\infty$ such that
$0 < n_m/(n_{m+1}-n_m) < c$ for some positive constant $c$.
Since
\[
\frac{n_m}{n_{m+1}} \pi'_m + \left(1- \frac{n_m}{n_{m+1}}\right) \delta_\theta
= \frac{n_m}{n_{m+1}} \pi_{n_m} + \left(1- \frac{n_m}{n_{m+1}}\right) \delta_\theta
\in \mathcal{P}_{\mu^q/n_{m+1}}^q
\]
for every $\theta \in \Theta^q$, where $\delta_\theta$ is the probability measure on $\Theta^q$ satisfying $\delta_\theta ( \{ \theta\} ) = 1$,
we have
\[
\tilde{\pi}_{m, \theta, u} := u \left\{ \frac{n_m}{n_{m+1}}\pi'_{m} + \left(1-\frac{n_m}{n_{m+1}}\right) \delta_\theta \right\}
+ (1- u) \pi'_{m+1}  \in \mathcal{P}_{\mu^q / n_{m+1}}^q
\]
for every $\theta \in \Theta^q$ and $0 \le u \le 1$.
Thus,
\begin{align*}
\frac{\partial}{\partial u} & D_q (\tilde{\pi}_{m,\theta,u}) \Bigg|_{u=0}
= \frac{\partial}{\partial u} \sum_{(x,y) \notin \mathcal{N}^q} p_{\tilde{\pi}_{m,\theta,u}} (x,y)
\log \frac{p_{\tilde{\pi}_{m,\theta,u}} (x,y)}{q(y;x) p_{\tilde{\pi}_{m,\theta,u}} (x)} \Bigg|_{u = 0} \notag \\
=& \sum_{(x,y) \notin \mathcal{N}^q} \biggl\{ \frac{\partial}{\partial u} p_{\tilde{\pi}_{m,\theta,u}} (x,y) \Bigg|_{u=0} \biggr\}
\log \frac{p_{\pi'_{m+1}}(x,y)}{q(y;x) p_{\pi'_{m+1}}(x)} \notag \\
=& \frac{n_m}{n_{m+1}} \sum_{(x,y) \notin \mathcal{N}^q} p_{\pi'_m} (x,y) \log \frac{p_{\pi'_{m+1}}(x,y)}{q(y;x) p_{\pi'_{m+1}}(x)}
- \sum_{(x,y) \notin \mathcal{N}^q} p_{\pi'_{m+1}} (x,y) \log \frac{p_{\pi'_{m+1}}(x,y)}{q(y;x) p_{\pi'_{m+1}}(x)} \notag \\
&+  \frac{n_{m+1}-n_m}{n_{m+1}}
\sum_{(x,y) \notin \mathcal{N}^q} p(x,y|\theta)
\log \frac{p_{\pi'_{m+1}}(x,y)}{q(y;x) p_{\pi'_{m+1}}(x)} \ge 0.
\end{align*}
Hence,
\begin{align*}
\lefteqn{
\sum_{(x,y) \notin \mathcal{N}^q} p(x,y|\theta)
\log \frac{p_{\pi'_{m+1}}(x,y)}{q(y;x) p_{\pi'_{m+1}}(x)}} \notag\\
\ge &\frac{n_{m+1}}{n_{m+1}-n_m} \sum_{(x,y) \notin \mathcal{N}^q} p_{\pi'_{m+1}} (x,y) \log \frac{p_{\pi'_{m+1}}(x,y)}{q(y;x) p_{\pi'_{m+1}}(x)}
- \frac{n_m}{n_{m+1}-n_m} \sum_{(x,y) \notin \mathcal{N}^q} p_{\pi'_m} (x,y) \log \frac{p_{\pi'_{m+1}}(x,y)}{q(y;x) p_{\pi'_{m+1}}(x)} \notag\\
=& \frac{n_{m+1}}{n_{m+1}-n_m} \sum_{(x,y) \notin \mathcal{N}^q} p_{\pi'_{m+1}} (x,y) \log \frac{p_{\pi'_{m+1}}(x,y)}{q(y;x) p_{\pi'_{m+1}}(x)} \notag\\
& + \frac{n_m}{n_{m+1}-n_m} \biggl\{ - \hspace{-10pt} \sum_{(x,y) \notin \mathcal{N}^q \cup \mathcal{N}^{\pi'_\infty}}
\hspace{-10pt} p_{\pi'_m} (x,y) \log \frac{p_{\pi'_{m+1}}(x,y)}{q(y;x) p_{\pi'_{m+1}}(x)}
- \hspace{-10pt} \sum_{(x,y) \in \mathcal{N}^{\pi'_\infty} \setminus \mathcal{N}^q} \hspace{-10pt} p_{\pi'_m} (x,y) \log p_{\pi'_{m+1}}(y|x) \notag\\
& + \hspace{-10pt} \sum_{(x,y) \in \mathcal{N}^{\pi'_\infty} \setminus \mathcal{N}^q} \hspace{-10pt} p_{\pi'_m} (x,y) \log q(y;x) \biggr\}
\end{align*}
\begin{align}
\ge& \frac{n_{m+1}}{n_{m+1}-n_m} \sum_{(x,y) \notin \mathcal{N}^q} p_{\pi'_{m+1}} (x,y) \log \frac{p_{\pi'_{m+1}}(x,y)}{q(y;x) p_{\pi'_{m+1}}(x)} \notag\\
& + \frac{n_m}{n_{m+1}-n_m} \biggl\{ - \hspace{-11pt} \sum_{(x,y) \notin \mathcal{N}^q \cup \mathcal{N}^{\pi'_\infty}} \hspace{-11pt}
p_{\pi'_m} (x,y) \log \frac{p_{\pi'_{m+1}}(x,y)}{q(y;x) p_{\pi'_{m+1}}(x)}
+ \hspace{-11pt} \sum_{(x,y) \in \mathcal{N}^{\pi'_\infty} \setminus \mathcal{N}^q} \hspace{-11pt} p_{\pi'_m} (x,y) \log q(y;x) \biggr\},
\label{eq:7-1}
\end{align}
where $\mathcal{N}^{\pi'_\infty} := \{(x,y) \in \mathcal{X} \times \mathcal{Y} \mid p_{\pi'_\infty}(x,y) = 0 \}$.
Here, we have
\begin{align}
\label{eq:lim1}
\lim_{m \rightarrow \infty} \hspace{-10pt} \sum_{(x,y)  \notin \mathcal{N}^q \cup \mathcal{N}^{\pi'_\infty}} \hspace{-10pt}
p_{\pi'_m} (x,y) \log \frac{p_{\pi'_{m+1}}(x,y)}{q(y;x) p_{\pi'_{m+1}}(x)}
= \hspace{-10pt} \sum_{(x,y)  \notin \mathcal{N}^q \cup \mathcal{N}^{\pi'_\infty}}
\hspace{-10pt} p_{\pi'_\infty} (x,y) \log \frac{p_{\pi'_\infty}(x,y)}{q(y;x) p_{\pi'_\infty}(x)},
\end{align}
because $p_{\pi'_\infty} (x,y) > 0$ for every $(x,y) \notin \mathcal{N}^{\pi'_\infty}$,
and
\begin{align}
\label{eq:lim2}
\lim_{m \rightarrow \infty} \hspace{-10pt} \sum_{(x,y) \in \mathcal{N}^{\pi'_\infty} \setminus \mathcal{N}^q} \hspace{-10pt} p_{\pi'_m} (x,y) \log q(y;x) = 0
= 
-  \hspace{-10pt} \sum_{(x,y) \in \mathcal{N}^{\pi'_\infty} \setminus \mathcal{N}^q} \hspace{-10pt} p_{\pi'_\infty} (x,y) \log \frac{p_{\pi'_\infty}(x,y)}{q(y;x) p_{\pi'_\infty}(x)}.
\end{align}
Therefore, from \eqref{eq:7-1}, \eqref{eq:lim1}, \eqref{eq:lim2}, and $0 < n_m/(n_{m+1}-n_m) < c$,
for every $\theta \in \Theta^q$,
\begin{align}
\liminf_{m \to \infty} & \sum_{(x,y) \notin \mathcal{N}^q} p(x,y|\theta) \log \frac{p_{\pi'_{m}}(x,y)}{q(y;x) p_{\pi'_{m}}(x)}
\ge \sum_{(x,y) \notin \mathcal{N}^q} p_{\pi'_\infty} (x,y) \log \frac{p_{\pi'_\infty}(x,y)}{q(y;x) p_{\pi'_\infty}(x)} \geq 0.
\label{keyineq}
\end{align}

By taking an appropriate subsequence $\{\pi''_k\}_{k=1}^\infty$ of $\{\pi'_m\}_{m=1}^\infty$, we can make
the sequences of real numbers $\{p_{\pi''_k}(y|x)\}^\infty_{k=1}$ converge for all $(x,y) \in \mathcal{X}^q \times \mathcal{Y}$
because $p_{\pi'_m}(x) > 0$ $(x \in \mathcal{X}^q)$ and $0 \leq p_{\pi'_m}(x,y)/p_{\pi'_m}(x) \leq 1$.

Then, from \eqref{keyineq}, if $\theta \in \Theta^q$,
\begin{align*}
R&(\theta, \lim\limits_{k \rightarrow \infty} p_{\pi_k''}(y|x)) = \sum_{x,y} p(x,y|\theta) \log \frac{p(y|x,\theta)}{\lim\limits_{k \rightarrow \infty} p_{\pi_k''}(y|x)}
= \sum_{(x,y) \notin \mathcal{N}^q} p(x,y|\theta) \log \frac{p(y|x,\theta)}{\lim\limits_{k \rightarrow \infty} p_{\pi_k''}(y|x)} \\
\leq& \sum_{(x,y) \notin \mathcal{N}^q} p(x,y|\theta) \log \frac{p(y|x,\theta)}{q(y;x)}
= \sum_{x,y} p(x,y|\theta) \log \frac{p(y|x,\theta)}{q(y;x)}
= R(\theta,q(y;x)) < \infty.
\end{align*}
Note that
the risk $R(\theta, \lim\limits_{k \rightarrow \infty} p_{\pi_k''}(y|x))$ does not depend on
the choice of $\lim\limits_{k \rightarrow \infty} p_{\pi''_k}(y|x)$ for $x \notin \mathcal{X}^q$,
although $\lim\limits_{k \rightarrow \infty} p_{\pi''_k}(y|x)$ is not uniquely determined for such $x$.

If $\theta \notin \Theta^q$, $R(\theta, q(y;x)) = \infty$ because
$- \hspace{-7pt} \sum\limits_{(x,y) \in \mathcal{N}^q} p(x,y|\theta) \log q(y;x) = \infty$.
For $x \notin \mathcal{X}^q$, $p(x|\theta) > 0$ only when $\theta \notin \Theta^q$.
Thus, if $x \notin \mathcal{X}^q$ is observed, then $R(\theta,q(y;x)) = \infty$
because $\theta \notin \Theta^q$.

Hence, the risk of the predictive density defined by
\begin{equation*}
\begin{cases}
\lim\limits_{k \rightarrow \infty} p_{\pi''_k} (y|x), & x \in \mathcal{X}^q \\
r (y;x), & x \notin \mathcal{X}^q,
\end{cases}
\end{equation*}
where $r(y;x)$ is an arbitrary predictive density, is not greater than that of $q(y;x)$ for every $\theta \in \Theta$.

Therefore, by taking a sequence $\{\varepsilon_n \in (0,1) \}_{n=1}^\infty$
that converges rapidly enough to $0$,
we can construct a predictive density
\begin{equation}
\lim_{k \to \infty} p_{\varepsilon_k \bar{\mu} + (1-\varepsilon_k)\pi''_k}(y|x) =
\begin{cases}
\lim\limits_{k \rightarrow \infty} p_{\pi''_k} (y|x), & x \in \mathcal{X}^q \\
p_{\bar{\mu}} (y|x), & x \notin \mathcal{X}^q
\end{cases}
\label{eq:1finalform}
\end{equation}
as a limit of Bayesian predictive densities based on priors
$\varepsilon_k \bar{\mu} + (1-\varepsilon_k) \pi''_k$,
where $\bar{\mu}$ is a measure on $\Theta$ such that $p_{\bar\mu} (x) > 0$ for every $x \in \mathcal{X}$.

Hence, the risk of the predictive density \eqref{eq:1finalform} is not greater than
that of $q(y;x)$ for every $\theta \in \Theta$.

\vspace{0.5cm}
\noindent
2)~ In this case, the proof becomes much simpler.
We assume that $\Theta^q \neq \emptyset$ because the assertion is obvious if $\Theta^q = \emptyset$.
Then, $D_q(\hat{\pi}) < \infty$ and $\hat{\pi}(\Theta^q) = 1$.
Thus, we can set $\mu^q = \hat{\pi}$ in the proof of 1).
Furthermore, we can set $\bar{\mu} = \hat{\pi}$ because
$p_{\hat{\pi}}(x) > 0$ for every $x \in \mathcal{X}$.
Therefore, the desired result can be proved without considering limits of Bayesian predictive densities.
\qed

\vspace{0.5cm}

We give two simple examples to clarify the meaning of Theorem 1 and its proof.

\vspace{0.5cm}

\noindent
Example 1.
Suppose that $\mathcal{X} = \{0, 1, 2\}$, $\mathcal{Y} = \{0,1\}$, $p(x,y|\theta) = {2 \choose x} \theta^x (1-\theta)^{2-x} \theta^y (1-\theta)^{1-y}$, and $\Theta = [0,1]$.
Let $q(y;x) =(x/2)^y (1-x/2)^{(1-y)}$, which is the plug-in predictive density with
the maximum likelihood estimate $\hat{\theta} = x/2$.
Then, $\mathcal{N}^q = \{(0, 1), (2, 0)\}$, $\Theta^q = \{0, 1\}$, and $\mathcal{X}^q = \{0,2\}$.
The prior defined by $\pi^{(w)} := w \delta_0 + (1-w) \delta_1 \in \mathcal{P}^q$ $(0 < w < 1)$ satisfies
\[
D_q(\pi^{(w)}) = \inf\limits_{\pi \in {\mathcal{P}^q}} D_q (\pi) = 0.
\]
We set $\mu^q = \pi^{(w)}$, which satisfies $p_{\mu^q}(x) > 0$ for $x \in \mathcal{X}^q$.
Then, we can set $\pi_n = \pi^{(w)}$ $(n=1,2,3,\ldots)$
because $\pi^{(w)} \in \mathcal{P}^q_{\mu^q/n}$ and $D_q(\pi^{(w)})=0$.
Then, $\lim\limits_{n \rightarrow \infty} p_{\pi_n}(y|x) = p_{\pi^{(w)}}(y|x)$.
Thus,
$\pi'_\infty = \pi^{(w)}$ and $\mathcal{N}^{\pi'_\infty} = \mathcal{N}^q$.

The prior $\pi^{(w)}$ does not specify the conditional density $p_{\pi^{(w)}}(y|x=1)$ because $p_{\pi^{(w)}}(x=1)=0$.
We set $\bar{\mu}(\dd \theta) = \dd \theta$ and
\[
\pi''_k = \frac{1}{k} \bar{\mu} + \biggl( 1 - \frac{1}{k} \biggr) \pi^{(w)}.
\]
Then, $\lim\limits_{k \rightarrow \infty} p_{\pi''_k}(y=0|x=0) = \lim\limits_{k \rightarrow \infty} p_{\pi''_k}(y=1|x=2) = 1$ and
$\lim\limits_{k \rightarrow \infty} p_{\pi''_k}(y=0|x=1) = \lim\limits_{k \rightarrow \infty} p_{\pi_k}(y=1|x=1) = 1/2$.
The risk function of the predictive density $\lim\limits_{k \rightarrow \infty} p_{\pi''_k}(y|x)$,
which is a limit of the Bayesian predictive densities,
is given by
\begin{align*}
R(\theta, \lim\limits_{k \rightarrow \infty} p_{\pi''_k}(y|x))
=
\begin{cases}
0, & \theta = 0 \in \Theta^q, \\
\infty, & \theta \in (0,1) = \Theta \setminus \Theta^q, \\
0, & \theta = 1 \in \Theta^q
\end{cases}
\end{align*}
and coincides with $R(\theta,q(y;x))$.
\qed

\vspace{0.5cm}
\noindent
Example 2.
Suppose that $\mathcal{X} = \{0,1,2\}$, $\mathcal{Y} = \{0,1\}$, $\Theta = \{ \theta_1, \theta_2\}$, 
$p((2,0)|\theta_1) = p((2,1)|\theta_1) = 0$, 
$p((0,0)|\theta_1) = p((1,1)|\theta_1) = 1/3$, 
$p((0,1)|\theta_1) = p((1,0)|\theta_1) = 1/6$,
$p((2,0)|\theta_2) = p((2,1)|\theta_2) = (1-\ep)/2$, and
$p((0,0)|\theta_2) = p((0,1)|\theta_2) = p((1,0)|\theta_2) = p((1,1)|\theta_2) = \ep/4$, where $0 < \ep < 1$. 

Consider a predictive density defined by $q(y=0;x=0) = q(y=1;x=1) = 2/3$,
$q(y=1;x=0) = q(y=0;x=1) = 1/3$, $q(y=0;x=2) = 1/3$, and $q(y=1;x=2) = 2/3$.
Then, $\mathcal{N}^q = \emptyset$, $\Theta^q = \Theta$, $\mathcal{P}^q = \mathcal{P}$, and $\mathcal{X}^q = \mathcal{X}$.

Then, $\hat{\pi} = \delta_{\theta_1}$ satisfies $D_q (\hat{\pi}) = \inf\limits_{\pi \in \mathcal{P}} D_q(\pi) = 0$
because $p(y|x,\theta_1) = q(y;x)$ except for the case $x=2$.
Since $p(x=2|\theta_1) = 0$, $p_{\hat{\pi}} (y|x =2)$ is not uniquely determined.
Thus, we consider a limit of Bayesian predictive densities.

Put $\mu = \delta_{\theta_1}/2 + \delta_{\theta_2}/2$.
It can be easily verified that $\pi_n = (1/n) \mu + (1-1/n) \delta_{\theta_1}$
satisfies $D_q(\pi_n) = \inf\limits_{\pi \in {\mathcal{P}_{\mu/n}}} D_q (\pi)$.
Then, $\lim\limits_{n \rightarrow \infty} p_{\pi_n} (y|x=0) = p(y|x=0,\theta_1) = q(y;x=0)$,
$\lim\limits_{n \rightarrow \infty} p_{\pi_n} (y|x=1) = p(y|x=1,\theta_1) = q(y;x=1)$,
$p_{\pi_n} (y|x=2) = p(y|x=2,\theta_2) \neq q(y;x=2)$.
By calculation,
we have $R(\theta_1, \lim\limits_{n \rightarrow \infty} p_{\pi_n}(y|x)) = R(\theta_1, q(y;x)) = 0$
and $R(\theta_2, \lim\limits_{n \rightarrow \infty} p_{\pi_n}(y|x)) = (\ep/2) \log (9/8) < R(\theta_2, q(y;x)) = (1/2) \log (9/8)$.
Thus, the performance of $\lim\limits_{n \rightarrow \infty} p_{\pi_n}(y|x)$
is better than that of $q(y;x)$ \qed
\section{\normalsize Latent information priors and minimax prediction}

In this section, we construct minimax predictive densities
that are limits of Bayesian predictive densities
based on prior sequences converging to latent information priors defined below.

A predictive density $q(y;x)$ is said to be minimax if it satisfies the equality
\[
\sup_{\theta \in \Theta} \sum_{x,y} p(x,y|\theta) \log \frac{p(y|x,\theta)}{q(y;x)} =
\inf_{\bar{q}} \sup_{\theta \in \Theta} \sum_{x,y} p(x,y|\theta) \log \frac{p(y|x,\theta)}{\bar{q}(y;x)}.
\]

The conditional mutual information between $y$ and $\theta$ given $x$ is defined by
\begin{align*}
I_{\theta,y |x} (\pi) := &\int \sum_{x,y} p(x,y|\theta) \log p(x,y|\theta) \dd\pi(\theta)  - \sum_{x,y} p_\pi (x,y) \log p_\pi(x,y) \\
&- \int \sum_x p(x|\theta) \log p(x|\theta)  \dd\pi( \theta) + \sum_x p_\pi (x) \log p_\pi (x),
\end{align*}
which is a function of $\pi \in \mathcal{P} $.
If $p_\pi(x) \neq 0$ for all $x \in \mathcal{X}$, then
\[
I_{\theta,y |x} (\pi) = \int \sum_{x,y} p(x,y|\theta) \log \frac{p(y|x,\theta)}{p_\pi (y|x)} \dd\pi( \theta).
\]
Since $u \log u \ \ (0 \le u \le 1)$ a bounded continuous function,
$I_{\theta,y |x} (\pi)$ is a bounded continuous function of $\pi \in \mathcal{P}$.

We define a latent information prior as a prior $\hat{\pi}$ that satisfies
$I_{\theta,y|x}(\hat{\pi}) = \sup\limits_{\pi \in \mathcal{P}} I_{\theta,y|x}(\pi)$.
Intuitively speaking, under the latent information prior,
the parameter $\theta$ has the maximum information about
the future observable $y$ under the condition that $x$ is observed.
Therefore, $\theta$ has the maximum amount of ``latent'' information, which we cannot observe through the data $x$.
Thus, the latent information prior corresponds to the ``worst case'' and is naturally related to minimaxity.
On the other hand, the minimum information prior discussed by Akaike (1983) is a prior maximizing the mutual information between the future observable $y$ and the data $x$.
This prior corresponds to the ``best case'' and is far from minimaxity.

The priors $\pi_\infty$ and $\hat{\pi}$ in Theorem 2 below are latent information priors.

\vspace{0.5cm}
\noindent
Theorem 2.
\begin{itemize}
\item[1)]
There exists a convergent prior sequence $\{\pi_n\}_{n=1}^\infty$
such that $\lim\limits_{n \to \infty} p_{\pi_n} (y|x)$ is a minimax predictive density
and the equality $I_{\theta,y|x}(\pi_\infty) = \sup\limits_{\pi \in \mathcal{P}} I_{\theta,y|x}(\pi)$
holds, where $\pi_\infty = \lim\limits_{n \to \infty} \pi_n$.
\item[2)]
Let $\hat{\pi} \in \mathcal{P}$ be a prior maximizing $I_{\theta,y|x}(\pi)$.
If $p_{\hat{\pi}} (x) > 0$ for all $x \in \mathcal{X}$, then $p_{\hat{\pi}}(y|x)$ is a minimax predictive density.
\qed
\end{itemize}

\vspace{0.5cm}
\noindent
Proof. 1)~
Let $\mu$ be a probability measure on $\Theta$ such that
$p_\mu (x) := \int p(x|\theta) \dd \mu(\theta) > 0$ for every $x \in \mathcal{X}$, and
let $\pi_n \in \mathcal{P}_{\mu/n} := \{\mu/n + (1-1/n)\pi \mid \pi \in \mathcal{P} \}$
be a prior satisfying
$I_{\theta,y|x} (\pi_n)= \sup\limits_{\pi \in \mathcal{P}_{\mu/n}}I_{\theta,y|x}(\pi)$.
From Lemma 2, there exists a convergent subsequence
$\{ \pi'_m \}^\infty_{m=1}$ of $\{ {\pi}_{n} \}^\infty_{n=1}$
such that $ I_{\theta,y|x}({\pi}'_\infty) = \sup\limits_{\pi \in \mathcal{P} } I_{\theta,y|x}(\pi) $,
where ${\pi}'_\infty = \lim\limits_{m \to \infty} \pi'_m$.
Let $n_m$ be the integer satisfying $\pi'_m = \pi_{n_m}$.
As in the proof of Theorem 1, we can take a subsequence $\{\pi'_m\}_{m=1}^\infty$ such that
$0 < n_m/(n_{m+1}-n_m) < c$ for some positive constant $c$.

Then, for every $\bar{\theta} \in \Theta$,
\[
\tilde\pi_{m, \btheta, u} := u \left\{\frac{n_m}{n_{m+1}} \pi'_m + \left( 1- \frac{n_m}{n_{m+1}} \right) \delta_{\bar{\theta}} \right\} + (1-u) \pi'_{m+1}
\]
belongs to $\mathcal{P}_{\mu/n_{m+1}}$ for $0 \le u \le 1$,
because $(n_m/n_{m+1}) \pi'_m + \left( 1 - n_m / n_{m+1} \right) \delta_\btheta \in \mathcal{P}_{\mu/n_{m+1}}$ and $\pi'_{m+1} \in \mathcal{P}_{\mu/n_{m+1}}$.

Thus,
\begin{align*}
\frac{\partial}{\partial u} & I_{\theta,y|x}  (\tilde\pi_{m, \btheta, u}) \Bigg|_{u = 0} =  \frac{\partial}{\partial u} 
\bigg( \int \sum_{x,y} p(x,y|\theta) \log p(x,y|\theta) \dd \tilde\pi_{m, \btheta, u} (\theta)
- \sum_{x,y} p_{\tilde\pi_{m, \btheta, u}}(x,y) \log p_{\tilde\pi_{m, \btheta, u}}(x,y) \\
& - \int \sum_x p(x|\theta) \log p(x|\theta) \dd \tilde\pi_{m, \btheta, u}(\theta)
+ \sum_x p_{\tilde\pi_{m, \btheta, u}} (x) \log p_{\tilde\pi_{m, \btheta, u}} (x) \bigg) \Bigg|_{u=0} \\
=& \frac{n_m}{n_{m+1}} \int \sum_{x,y} p(x,y|\theta) \log p(x,y|\theta) \dd \pi'_{m} (\theta)
+ \biggl(1 - \frac{n_m}{n_{m+1}}\biggr) \sum_{x,y} p(x,y|\bar{\theta}) \log p(x,y|\bar{\theta}) \\
& - \int \sum_{x,y} p(x,y | \theta) \log p(x,y|\theta) \dd \pi'_{m+1} (\theta)
- \sum_{x,y}  \frac{\partial}{\partial u} p_{\tilde{\pi}_{m,\btheta,u}}(x,y) \Bigg|_{u=0}  \log p_{\pi'_{m+1}} (x,y) \\
& - \frac{n_m}{n_{m+1}} \int \sum_x p(x|\theta) \log p(x|\theta) \dd \pi'_{m} (\theta)
- \biggl(1 - \frac{n_m}{n_{m+1}} \biggr) \sum_x p(x|\bar{\theta}) \log p(x|\bar{\theta}) \\
& + \int \sum_x p(x|\theta) \log p(x|\theta) \dd \pi'_{m+1} (\theta)
+ \sum_x  \frac{\partial}{\partial u} p_{\tilde{\pi}_{m,\btheta,u}}(x) \Bigg|_{u=0}  \log p_{\pi'_{m+1}} (x) \\
= & \biggl(1 - \frac{n_m}{n_{m+1}} \biggr) \sum_{x,y} p(x,y|\bar{\theta}) \log \frac{p(x,y|\bar{\theta})}{p(x|\bar{\theta})}
- \biggl(1 - \frac{n_m}{n_{m+1}} \biggr) \sum_{x,y} p(x,y|\bar{\theta}) \log \frac{p_{\pi'_{m+1}} (x,y)}{p_{\pi'_{m+1}} (x)} \\
& + \frac{n_m}{n_{m+1}} \int \sum_{x,y} p(x,y|\theta) \log \frac{p(x,y|\theta)}{p(x|\theta)} \dd \pi'_{m} (\theta)
- \int \sum_{x,y} p(x,y | \theta) \log \frac{p(x,y|\theta)}{p(x|\theta)} \dd \pi'_{m+1} (\theta)\\
&- \frac{n_m}{n_{m+1}} \sum_{x,y} p_{\pi'_{m}} (x,y) \log \frac{p_{\pi'_{m+1}} (x,y)}{p_{\pi'_{m+1}} (x)}
+ \sum_{x,y} p_{\pi'_{m+1} }(x,y) \log \frac{p_{\pi'_{m+1}} (x,y)}{p_{\pi'_{m+1}} (x)} \le 0 ,
\end{align*}
where we used
\[
\frac{\partial}{\partial u} p_{\tilde\pi_{m, \btheta, u}}(x,y) = \frac{n_m}{n_{m+1}} p_{\pi'_m} (x,y)
+ \biggl( 1-\frac{n_m}{n_{m+1}} \biggr) p(x,y|\bar\theta) - p_{\pi'_{m+1}}(x,y).
\]

Noting that $p_{\pi'_m}(x) > 0$ for every $m$ and $x \in \mathcal{X}$ and that
$p(x,y|\theta) \log p(y|x,\theta)=0$ if $p(x|\theta) = 0$, we have
\begin{align*}
& \biggl(1 - \frac{n_m}{n_{m+1}}\biggr) \sum_{x,y} p(x,y|\bar{\theta})
\log \frac{p(y|x,\bar{\theta})}{p_{\pi'_{m+1}} (y|x)}
+ \frac{n_m}{n_{m+1}} \int \sum_{x,y} p(x,y|\theta)
\log \frac{p(y|x,\theta)}{p_{\pi'_{m+1}} (y|x)} \dd \pi'_{m} (\theta) \\
& - \int \sum_{x,y} p(x,y|\theta)
\log\frac{p(y|x,\theta)}{ p_{\pi'_{m+1}} (y|x)} \dd \pi'_{m+1} (\theta) \le 0.
\end{align*}
Hence,
\begin{align}
\sum_{x,y} & p(x,y| \bar{\theta}) \log \frac{p(y|x, \bar{\theta})}{p_{\pi'_{m+1}}(y|x)}
\leq
- \frac{n_m}{n_{m+1}-n_m} \biggl\{ \int \sum_{(x,y) \notin \mathcal{N}^{\pi'_\infty}} p(x,y|\theta) \log \frac{p(y|x,\theta)}{p_{\pi'_{m+1}} (y|x)} \dd \pi'_{m} (\theta) \notag\\
& + \int \sum_{(x,y) \in \mathcal{N}^{\pi'_\infty}} p(x,y|\theta) \log p(y|x,\theta) \dd \pi'_{m} (\theta)
- \int \sum_{(x,y) \in \mathcal{N}^{\pi'_\infty}} p(x,y|\theta) \log p_{\pi'_{m+1}} (y|x) \dd \pi'_{m} (\theta) \biggr\} \notag\\
& + \frac{n_{m+1}}{n_{m+1}-n_m} \int \sum_{x,y} p(x,y|\theta) \log\frac{p(y|x,\theta)}{ p_{\pi'_{m+1}} (y|x)} \dd \pi'_{m+1} (\theta) \notag\\
\leq
& - \frac{n_m}{n_{m+1}-n_m} \biggl\{ \int \hspace{-10pt} \sum_{(x,y) \notin \mathcal{N}^{\pi'_\infty}} \hspace{-10pt}
p(x,y|\theta) \log \frac{p(y|x,\theta)}{p_{\pi'_{m+1}} (y|x)} \dd \pi'_{m} (\theta)
+ \int \hspace{-10pt} \sum_{(x,y) \in \mathcal{N}^{\pi'_\infty}} \hspace{-10pt} p(x,y|\theta) \log p(y|x,\theta) \dd \pi'_{m} (\theta) \biggr\} \notag\\ &
+ \frac{n_{m+1}}{n_{m+1}-n_m} \int \sum_{x,y} p(x,y|\theta) \log\frac{p(y|x,\theta)}{ p_{\pi'_{m+1}} (y|x)} \dd \pi'_{m+1} (\theta), \label{A}
\end{align}
where $\mathcal{N}^{\pi'_\infty} := \{(x,y) \in \mathcal{X} \times \mathcal{Y} \mid p_{\pi'_\infty}(x,y) = 0 \}$.
Here, we have
\begin{align}
\lim_{m \to \infty}
\int \hspace{-10pt} \sum_{(x,y) \notin \mathcal{N}^{\pi'_\infty}} \hspace{-10pt}
p(x,y|\theta) \log \frac{p(y|x,\theta)}{p_{\pi'_{m+1}} (y|x)} \dd \pi'_{m} (\theta)
= \int \hspace{-10pt} \sum_{(x,y) \notin \mathcal{N}^{\pi'_\infty}} \hspace{-10pt}
p(x,y|\theta) \log \frac{p(y|x,\theta)p_{\pi'_\infty}(x)}
{p_{\pi'_\infty}(x,y)} \dd \pi'_{\infty} (\theta)
\label{B}
\end{align}
and
\begin{align}
\lim_{m \to \infty}
\int \hspace{-10pt} \sum_{(x,y) \in \mathcal{N}^{\pi'_\infty}} \hspace{-10pt} p(x,y|\theta) \log p(y|x,\theta) \dd \pi'_{m} (\theta)
&= \int \hspace{-10pt} \sum_{(x,y) \in \mathcal{N}^{\pi'_\infty}} \hspace{-10pt} p(x,y|\theta) \log p(y|x,\theta) \dd \pi'_{\infty} (\theta) \notag\\
= \int \hspace{-10pt} & \sum_{(x,y) \in \mathcal{N}^{\pi'_\infty}} \hspace{-10pt} p(x,y|\theta)
\log \frac{p(y|x,\theta)p_{\pi'_\infty}(x)}{p_{\pi'_{\infty}}(x,y)}
\dd \pi'_{\infty} (\theta) = 0,
\label{C}
\end{align}
because $p(x,y|\theta) \log p(x,y|\theta)$ and $p(x|\theta) \log p(x|\theta)$
are bounded continuous functions of $\theta$ for every fixed $(x,y)$.

From \eqref{A}, \eqref{B}, \eqref{C}, and $0 < n_m/(n_{m+1}-n_m) < c$, we have, for every
$\bar{\theta} \in \Theta$,
\[
\limsup_{m \to \infty} \sum_{x,y} p(x,y| \bar{\theta}) \log \frac{p(y|x, \bar{\theta})}{p_{\pi'_m}(y|x)} \le
\int \sum_{x,y} p(x,y|\theta) \log 
\frac{p(y|x, \theta)p_{\pi'_\infty}(x)}{p_{\pi'_\infty}(x,y)} \dd \pi'_\infty(\theta).
\]
By taking an appropriate subsequence $\{\pi''_k\}_{k=1}^\infty$ of $\{\pi'_m\}_{m=1}^\infty$, 
we can make $\{p_{\pi''_k} (y|x)\}_{k=1}^\infty$ converges for every $(x,y)$ as $k \to \infty$. 
Then,
for every $\bar{\theta} \in \Theta$,
\begin{align}
\sum_{x,y} & p(x,y |\bar{\theta}) \log \frac{p(y|x,\bar{\theta})} {\displaystyle \lim_{k \to \infty} p_{\pi''_k}(y|x)}
\le \int \sum_{x,y} p(x,y|\theta) \log \frac{p(y|x,\theta)}{\displaystyle\lim_{k \to \infty} p_{\pi''_k}(y|x)} \dd \pi''_\infty (\theta),
\label{jun}
\end{align}
where $\pi''_\infty = \pi'_\infty = \lim\limits_{k \rightarrow \infty} \pi''_k$,
because $\lim\limits_{k \rightarrow \infty} p_{\pi''_k}(y|x) = p_{\pi''_\infty}(y|x)$
for $x$ with $p_{\pi''_\infty}(x) > 0$.

On the other hand, we have
\begin{align}
\int \sum_{x,y} & p(x,y|\theta) \log \frac{p(y|x,\theta)}{\displaystyle\lim_{k \to \infty} p_{\pi''_k}(y|x)} \dd \pi''_\infty(\theta)
= \inf_q \int \sum_{x,y} p(x,y|\theta) \log \frac{p(y|x,\theta)}{q(y;x)} \dd \pi''_\infty (\theta) \notag \\
\leq& \sup_{\pi \in \mathcal{P}} \inf_q \int \sum_{x,y} p(x,y|\theta) \log \frac{p(y|x,\theta)}{q(y;x)} \dd \pi (\theta)
\leq \inf_q \sup_{\pi \in \mathcal{P}} \int \sum_{x,y} p(x,y|\theta) \log \frac{p(y|x,\theta)}{q(y;x)} \dd \pi (\theta) \notag \\
=& \inf_q \sup_{\theta \in \Theta} \sum_{x,y} p(x,y|\theta) \log \frac{p(y|x,\theta)}{q(y;x)}
\leq \sup_{\theta \in \Theta} \sum_{x,y} p(x,y|\theta) \log \frac{p(y|x,\theta)}{\displaystyle\lim_{k \to \infty} p_{\pi''_k}(y|x)}.
\label{gyaku}
\end{align}
The first equality is because the Bayes risk
\[
\int R(\theta;q(y;x)) \pi''_\infty( \dd \theta ) =
\int
\sum_{x,y} p(x,y|\theta) \log \frac{p(y|x,\theta)}{q(y;x)} \pi''_\infty(\dd \theta)
\]
with respect to $\pi''_\infty \in \mathcal{P}$ is minimized when
\[
q(y;x) = p_{\pi''_\infty}(y|x)
:= \frac{\int p(x,y|\theta) \pi''_\infty(\dd \theta)}{\int p(x|\theta) \pi''_\infty(\dd \theta)};
\]
see Aitchison (1975).
Although $p_{\pi''_\infty}(y|x)$ is not uniquely determined for $x$ with $p_{\pi''_\infty}(x) = 0$,
the Bayes risk does not depend on the choice of $p_{\pi''_\infty}(y|x)$ for such $x$.

From \eqref{jun} and \eqref{gyaku}, we have
\[
\inf_q \sup_{\theta \in \Theta} \sum_{x,y} p(x,y|\theta) \log \frac{p(y|x,\theta)}{q(y;x)}
= \sup_{\theta \in \Theta} \sum_{x,y} p(x,y|\theta) \log \frac{p(y|x,\theta)}{\displaystyle\lim_{k \to \infty} p_{\pi''_k}(y|x)}.
\]
Therefore, the predictive density $\lim\limits_{k \rightarrow \infty} p_{\pi''_k}(y|x)$ is minimax.

\vspace{0.5cm}
\noindent
2)~ In this case, the proof becomes much simpler.
By setting $\mu = \hat{\pi}$ in the proof of 1), we have $\pi_n = \hat{\pi}$ $(n=1,2,3,\ldots)$.
Thus, $\lim\limits_{n \rightarrow \infty} p_{\pi_n}(y|x) = p_{\hat{\pi}}(y|x)$,
and the desired result can be proved without considering limits of Bayesian predictive densities.
\qed

\section {\normalsize Numerical results and discussions}

Let $p(x|\theta) = {N \choose x} \theta^x (1-\theta)^{N-x}$ $(x=0,1,\ldots,N)$,
$p(y|\theta) = {M \choose y} \theta^y (1-\theta)^{M-y}$ $(y=0,1,\ldots,M)$,
and $\Theta = \{ 0.1 k \mid k = 0,1,2,\ldots,10 \}$
in which $\theta$ takes a value.
Although this example is relatively simple in the sense that $x$ and $y$ are independent given $\theta$,
the behavior of priors is not trivial.

The latent information priors, which maximize $I_{\theta,y|x}(\pi)$,
for 16 sets of values of $(N,M)$ are obtained numerically; see Figure 1.
\begin{figure}[h]
\begin{center}
\includegraphics[scale=0.55]{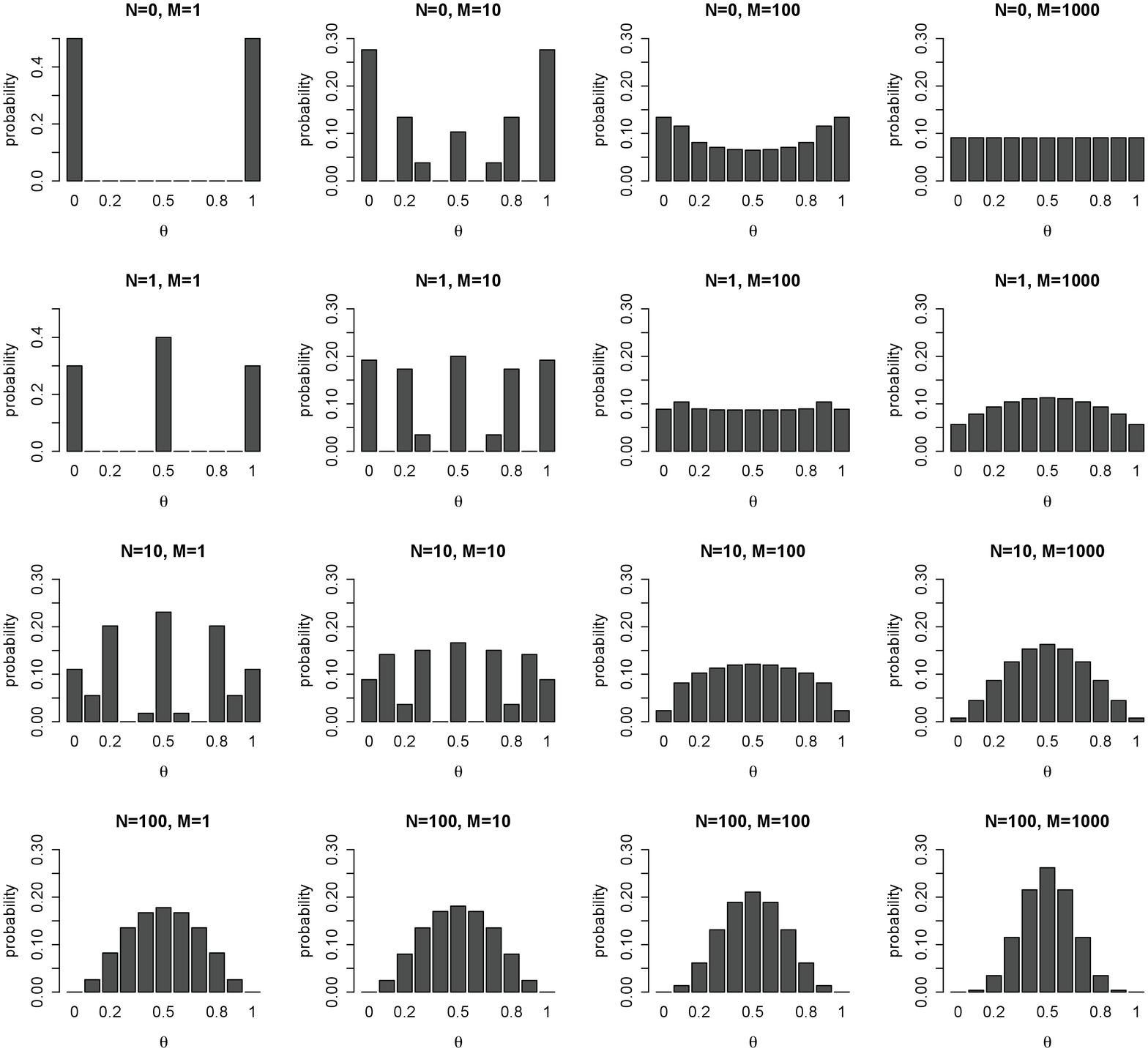}
\label{prior}

\noindent
{Figure 1. Latent information priors for various $(N,M)$ values}
\end{center}
\end{figure}

The prior for $(N,M)=(0,1000)$ is almost uniform and is similar to the reference prior
because the reference prior is the latent information prior with $N=0$ and $M \rightarrow \infty$.
It is widely known the reference prior is uniform when the parameter space is a finite set.
The latent information prior for $(N,M) = (0,100)$ is similar to the histogram of
the Jeffreys prior density $\theta^{-1/2}(1-\theta)^{-1/2}/B(1/2,1/2)$
for the binomial model with the ordinary parameter space $\Theta = [0,1]$.

When both of $N$ and $M$ are small the priors assign weights
only on a limited number of points in $\Theta$.
This corresponds to the phenomenon concerning the $k$-reference prior
studied by Berger, Bernardo, and Mendoza~(1989).
The $k$-reference prior is the latent information prior with $N=0$ and $M=k$.

When $N$ is large, the priors assign more weights to parameter values close to $0.5$.
The shapes of priors are quite different from the uniform density or the histogram of the Jeffreys prior for the binomial model
with the ordinary parameter space $\Theta = [0,1]$.

These observations show that the latent information priors strongly depend on $(N,M)$.
This indicates that we need to abandon the context invariance (see Dawid (1983)) of priors.

The relation between the conditional mutual information and predictive densities parallels to
that between the unconditional mutual information and Bayes codes in information theory
except for the care for the case $ p_\pi (x) =0$.
Many studies on the unconditional mutual information and minimax prediction and coding
have been carried out; see, for example,
Ibragimov and Hasminskii (1972), Gallager (1979),
Davisson and Leon-Garcia (1980), Clark and Barron (1994), and Haussler (1997).
See also Gr\"{u}nwald and Dawid (2004) for discussions in a very general setting.
The conditional mutual information $I_{\theta, y|x}(\pi)$ coincides with
the Bayes risk of the Bayesian predictive density based on $\pi$.
Therefore, it is natural that the prior maximizing $I_{\theta, y|x} (\pi)$
corresponds to minimax prediction based on data.

In general, the priors based on the unconditional mutual information and that based on the conditional mutual information are quite different.
Latent information priors maximizing the conditional mutual information
could play important roles in statistical applications.
Although we have discussed submodels of multinomial models,
essential part of our discussion seem to hold for more
general models under suitable regularity conditions including compactness of the model
as in the theory based on the unconditional mutual information studied by Haussler (1997).

The explicit forms of latent information priors are usually complex and difficult to obtain
unless the parameter space is finite.
For actual applications,
it is important to develop approximation methods and asymptotic theory
in various settings other than the situation $N=0, M \rightarrow \infty$
studied in the reference analysis.
When $I_{\theta, y|x}(\pi)$ is close to $I_{\theta, y|x}(\hat{\pi})$, a prior $\pi$ is considered to be close to $\hat{\pi}$ 
because $I_{\theta, y|x}(\pi)$ is a concave function of $\pi$.
These topics require further research and will be discussed in other places.

\vspace{0.5cm}

\begin{center}
References
\end{center}

{\small

\refset
Aitchison, J.~(1975).~Goodness of prediction fit. {\it Biometrika}, {\bf 62}, 547--554. 

\refset
Akaike, H.~(1983).~On minimum information prior distributions.
{\it Annals of the Institute of Statistical Mathematics}, {\bf 35} Part A, 139--149.

\refset
Berger, J.~O., Bernardo, J.~M., and Mendoza, M.~(1989). On priors that maximize expected information.
{\it Recent Developments of Statistics and its Applications}, Klein, J.~and Lee, J.~eds.,
Freedom Academy, Seoul 1--20. 

\refset
Bernardo,~J.~M.~(1979).~Reference posterior distributions for
Bayesian inference (with discussion).\
{\it Journal of Royal Statistical Society} B, {\bf 41},~113--147.

\refset
Bernardo,~J.~M.~(2005).~Reference analysis.
{\it Handbook of Statistics}, {\bf 25}, 
Dey, K.\ K.\ and Rao C.\ R.\ eds., Elsevier, Amsterdam 17--90.

\refset
Clarke, B.~S.~and Barron,~A.~R.~(1994).~Jeffreys' prior is asymptotically least favorable under
entropy risk.\ {\it Journal of Statistical Planning and Inference}, {\bf 41},~36--60.

\refset
Davisson, L.~and Leon-Garcia, A.~(1980).~A source matching approach to finding minimax codes.\
{\it IEEE Transactions on Information Theory}, {\bf 26}, 166--174.

\refset
Dawid, A.~P.~(1983).~Invariant Prior Distributions. {\it Encyclopedia of Statistical Sciences}, {\bf 4},
Kotz, S.~Johnson,  N.~L.~and Read C.~B.~eds., Wiley-Interscience, New York 228--236.

\refset
Gallager, R.~(1979).~Source coding with side information and universal coding.\
Technical Report LIDSP-937, M.I.T. Laboratory for Information and Decision Systems.

\refset
Geisser, S.~(1979).~Discussion on ``Reference posterior distributions for Bayesian inference''
by J.~M.~Bernardo, {\it Journal of Royal Statistical Society} B, {\bf 41},~136--137.

\refset
Gr\"{u}nwald, P.~D.~and Dawid, A.~P.~(2004).~\
Game theory, maximum entropy, minimum discrepancy and robust {B}ayesian decision theory.\
{\it Annals of Statistics}, {\bf 32},~1367--1433.

\refset
Haussler, D.~(1997).~A general minimax result for relative entropy.\
{\it IEEE Transactions on Information Theory}, {\bf 43}, 1276--1280.

\refset
Ibragimov, I.~A.~and Hasminskii, R.~Z.~(1973).~On the information contained in a sample about a parameter.\
In {\it 2nd Intl. Symp. on Information Theory}, Akademiai, Kiado, Budapest 295--309.

\refset
Komaki, F. ~(2004).~Simultaneous prediction of independent Poisson observables.\
{\it Annals of Statistics}, {\bf 32}, 1744--1769.

\refset
Kuboki, H. ~(1998).~Reference priors for prediction.\
{\it Journal of Statistical Planning and Inference}, {\bf 69}, 295--317.
}
\end{document}